\documentclass[12pt]{amsart}

\usepackage{amsmath, amscd, graphicx, latexsym, hyperref, rlepsf, times}

\oddsidemargin 0in \theoremstyle{plain} \topmargin 0in

\oddsidemargin .25in \theoremstyle{plain}

\newtheorem{Thm}{Theorem}[section]
\newtheorem{Lem}[Thm]{Lemma}
\newtheorem{Cor}[Thm]{Corollary}
\newtheorem{Prop}[Thm]{Proposition}
\newtheorem{Def}[Thm]{Definition}
\newtheorem{Rem}[Thm]{Remark}

\newcommand{\pob}{partial open book decomposition }

\newcommand{\bfz}{{\mathbb{Z}}}
\newcommand{\bfc}{{\mathbb{C}}}

\def\p{\partial}

\def\v{\vskip.12in}
\def\a{\alpha}
\def\A{$\boldmath$\alpha$\unboldmath$}
\def\B{$\boldmath$\beta$\unboldmath$}
\def\Atek{\boldmath$\alpha \ $\unboldmath}
\def\Btek{\boldmath$\beta \ $\unboldmath}
\def\b{\beta}
\def\d{\delta}
\def\G{\Gamma}
\def\g{c}

\def\S{\Sigma}
\def\eps{\epsilon}
\def\SF{SFH(-M,-\G)}

\begin{document}

\title[Partial open books and the contact class in sutured Floer homology]{Partial
open book decompositions \\ and the contact class in sutured Floer homology}

\author{Tolga Etg\"u}

\author{Burak Ozbagci}

\begin{abstract}

We demonstrate how to combinatorially calculate the EH-class of a
compatible contact structure in the sutured Floer homology group
of a balanced sutured three manifold which is associated to an
\emph{abstract} partial open book decomposition. As an application
we show that every contact three manifold (closed or with convex
boundary) can be obtained by gluing tight contact handlebodies
whose EH-classes are nontrivial.

\end{abstract}

\address{Department of Mathematics \\ Ko\c{c} University \\ Istanbul, Turkey}
\email{tetgu@ku.edu.tr} \email{bozbagci@ku.edu.tr}
\subjclass[2000]{}


\thanks{}

\v \v \v

\maketitle

\setcounter{section}{-1}


\section{Introduction}

A sutured manifold $(M,\G)$ is a compact oriented $3$-manifold
with nonempty boundary, together with a compact subsurface $\G=
A(\G) \cup T(\G) \subset \p M$, where $A(\G)$ is a union of
pairwise disjoint annuli and $T(\G)$ is a union of tori. Moreover
each component of $\p M \setminus \G$ is oriented, subject to the
condition that the orientation changes every time we nontrivially
cross $A(\G)$. Let $R_+ (\G)$ (resp. $R_- (\G)$) be the open
subsurface of $\p M \setminus \G$ on which the orientation agrees
with (resp. is the opposite of ) the boundary orientation on $\p
M$. A sutured manifold $(M,\G)$ is balanced if $M$ has no closed
components, $\pi_0(A(\G)) \to \pi_0(\p M) $ is surjective, and
$\chi(R_+ (\G))= \chi (R_- (\G))$ on every component of $M$. It
follows that if $(M,\G)$ is balanced, then $\G=A(\G)$ and every
component of $\p M$ nontrivially intersects $\G$. Since all our
sutured manifolds will be balanced in this paper, we can think of
$\G$ as a set of \emph{oriented curves} on $\p M$ by identifying
each annulus in $\G$ with its core circle. Here $\Gamma$ is
oriented as the boundary of $R_+ (\G)$.

Let $\xi$ be a contact structure on a compact oriented
$3$-manifold $M$ whose dividing set on the convex boundary $\p M$
is denoted by $\G$. Then it is not too hard to see that $(M,
\Gamma)$ is a \emph{balanced} sutured manifold with the
identification we mentioned above. Recently, Honda, Kazez and
Mati\'{c} \cite{hkm1} introduced an invariant $EH(M,\G,\xi)$ of
the contact structure $\xi$ which lives in the sutured Floer
homology group $\SF$ defined for the balanced sutured manifold
$(M, \Gamma)$ by Juh\'{a}sz \cite{juh}. This invariant generalizes
the contact class in Heegaard Floer homology in the closed case as
defined by Ozsv\'{a}th and Szab\'{o} \cite{os} and reformulated in
\cite{hkm}.

In order to define $EH(M,\G,\xi)$, Honda, Kazez and Mati\'{c}
first construct a \pob of $M$ compatible with the given contact
structure $\xi$ by generalizing the work of Giroux \cite{g} in the
closed case. Then they obtain an admissible balanced  Heegaard
diagram for $(-M,-\G)$ which not only leads to the calculation of
the sutured Floer homology group $\SF$ but also includes the
description (similar to the one in the closed case again due to
Honda, Kazez and Mati\'{c} \cite{hkm}) of a certain cycle
descending to the contact class $EH (M, \G, \xi)$ in $\SF$, in
fact in $\SF/\{\pm 1\}$, but this $\pm 1$ ambiguity is usually
suppressed.

On the other hand, in \cite{eo} the authors gave an abstract
definition of a partial open book decomposition of a compact
$3$-manifold with boundary; associated a balanced sutured manifold
to a partial open book decomposition and constructed a compatible
contact structure on this sutured manifold whose dividing set on
the convex boundary agrees with the suture. In this paper we show that the
sutured Floer homology group and the EH-contact class can be
combinatorially calculated starting from an \emph{abstract}
partial open book decomposition.

We include here several sample calculations of the EH-contact
class and as an application we show that every contact three
manifold (closed or with convex boundary) can be obtained by
gluing tight contact handlebodies whose EH-classes are nontrivial.

The reader is advised to turn to \cite{hkm1} and \cite{eo} for
necessary background on partial open book decompositions, to
Juh\'{a}sz's papers \cite{juh} and \cite{juh1} for the definition
and properties of the sutured Floer homology of balanced sutured
manifolds and to Etnyre's notes \cite{e} for the related material
on contact topology of three-manifolds.

\section{Partial open book decompositions and compatible contact structures}\label{def}

In this section we quickly review basics about partial open book
decompositions and compatible contact structures as described in
\cite{hkm1} and  \cite{eo}.

\begin{Def} [\cite{eo}] \label{Pob} A \pob is a triple $(S,P,h)$ satisfying the
following conditions:

$(1)$ $S$ is a compact oriented connected surface with $\p S \neq
\emptyset$,

$(2)$ $P= P_1 \cup P_2 \cup \ldots \cup P_r$ is a proper (not
necessarily connected) subsurface of $S$ such that $S$ is obtained
from $\overline{S \setminus P}$ by successively attaching $1$-handles $P_1,
P_2, \ldots, P_r$,

$(3)$ $h:P \to S$ is an embedding such that $h|_A =$ identity, where
$A=\p P\cap \p S$.

\end{Def}

Given a \pob $(S,P,h)$, we construct a sutured manifold $(M,\G)$
as follows: Let $$H=(S \times [-1,0])/\sim
$$ where $(x,t) \sim (x,t')$ for $x \in \p S$ and $t, t' \in
[-1,0]$. It is easy to see that $H$ is a solid handlebody whose
oriented  boundary is the surface $S \times \{0 \} \cup - S \times
\{-1\}$ (modulo the relation $(x,0) \sim (x,-1)$ for every
$x \in \p S$). Similarly let $$N=(P \times [0,1])/\sim$$ where $(x,t)
\sim (x,t')$ for $x \in A$ and $t, t' \in [0,1]$. Observe that each component of
$N$ is also a solid handlebody. The oriented boundary of $N$ can be
described as follows: Let the arcs $\g_1, \g_2, \ldots, \g_n$
denote the connected components of $\overline{\p P \setminus \p
S}$. Then, for $1 \leq i \leq n$, the disk $D_i= (\g_i \times
[0,1])/\sim$ belongs to $\p N$. Thus part of $\p N$ is given by
the disjoint union of $D_i$'s. The rest of $\p N$ is the surface
$P \times \{1 \} \cup - P \times \{0\}$ (modulo the relation
$(x,0) \sim (x,1)$ for every $x \in A$).

Let $M= N \cup H$ where we glue these manifolds by identifying
$P\times \{0\} \subset \p N $ with $P \times \{0\} \subset \p H $
and $P \times \{1\} \subset \p N $ with $h(P) \times \{-1\}
\subset \p H $. Since the gluing identification is orientation
reversing $M$ is a compact oriented $3$-manifold with oriented boundary
$$\p M =
(S\setminus P ) \times \{ 0\}
\cup -(S \setminus h(P)) \times \{ -1\}
\cup (\overline{\p P \setminus \p S}) \times [0,1]
$$
(modulo the identifications given above).

We define the suture $\G$ on $\p M$ as the set of closed curves
obtained by gluing the arcs $\g_i \times \{1/2 \} \subset \p N$,
for $1 \leq i \leq n$, with the arcs in $(\overline{\p S \setminus
\p P}) \times \{0\} \subset \p H$, hence as an oriented simple
closed curve and modulo identifications
$$ \G =
(\overline{\p S \setminus \p P}) \times \{0\}
\cup
- (\overline{\p P \setminus \p S}) \times \{1/2\}
\ .$$

In \cite{eo} we showed that the sutured manifold $(M,\G)$
associated to a \pob $(S,P,h)$ is balanced.

{\Prop [\cite{eo}] \label{torisu} Let $(M,\G)$ be the balanced
sutured manifold associated to a \pob $(S,P,h)$. Then  there
exists a contact structure $\xi$ on $M$ satisfying the following
conditions:

$(1)$ $\xi$ is tight when restricted to $H$ and $N$,

$(2)$ $\p H$ is a convex surface in $(M, \xi)$ whose dividing set
is $\p S \times \{0 \}$,

$(3)$ $\p N $ is a convex surface in $(M, \xi)$ whose dividing set
is $\p P \times \{ 1/2 \}$.

Moreover such $\xi$ is unique up to isotopy.}

\vspace{0.2in}

Let $(M,\G)$ be the balanced sutured manifold associated to a \pob
$(S,P,h)$. A contact structure $\xi$ on $(M,\G)$ is said to be
compatible with $(S,P,h)$ if it satisfies conditions $(1),(2)$ and
$(3)$ stated in Proposition~\ref{torisu}.  It follows from
Proposition~\ref{torisu} that every \pob has a unique compatible
contact structure, up to isotopy, on the balanced sutured manifold
associated to it, such that the dividing set of the convex
boundary is isotopic to the suture.

Two partial open book decompositions $(S,P,h)$ and
$(\widetilde{S}, \widetilde{P}, \widetilde{h})$ are isomorphic if
there is a diffeomorphism $f: S \to \widetilde{S}$ such that $f(P)
= \widetilde{P}$ and $\widetilde{h} = f \circ h \circ
(f^{-1})|_{\widetilde{P}}$. Consequently, if $(S,P,h)$ and
$(\widetilde{S}, \widetilde{P}, \widetilde{h})$ are isomorphic
partial open book decompositions, then the associated compatible
contact $3$-manifolds $(M, \G, \xi)$ and $(\widetilde{M},
\widetilde{\G}, \widetilde{\xi})$ are isomorphic.

The following theorem is the key to
obtaining a description of a partial open book decomposition of
$(M,\G, \xi)$ in the sense of Honda, Kazez and Mati\'{c}.

{\Thm [\cite{hkm1}, Theorem 1.1]  \label{honkazmat} Let $(M,\G)$
be a balanced sutured manifold and let $\xi$ be a contact
structure on $M$ with convex boundary whose dividing set $\G_{\p
M}$ on $\p M$ is isotopic to $\G$. Then there exist a Legendrian
graph $K \subset M $ whose endpoints lie on $\G \subset \p M$ and
a regular neighborhood $N(K) \subset M $ of $K$ which satisfy the
following:

\begin{itemize}
\item[](A) \ (i) $T=\overline{\p N(K) \setminus \p M}$ is a convex surface
with Legendrian boundary.
\begin{itemize}
\item[](ii) For each component $\gamma_i$ of
$\p T$, $\gamma_i \cap \G_{\p M}$ has two connected components.
\item[](iii) There is a system of pairwise disjoint compressing disks
$D^\a_j$ for $N(K)$ so that $\p D^\a_j$ is a curve on $T$
intersecting the dividing set $\G_T$ of $T$ at two points and each
component of $N(K) \setminus \cup_j D^\a_j$ is  a standard contact
$3$-ball, after rounding the edges.
\end{itemize}
\item[](B) \ (i) Each component $H$ of $\overline{M \setminus N(K)}$ is a
handlebody (with convex boundary).
\begin{itemize}
\item[](ii) There is a system of
pairwise disjoint compressing disks $D^{\d}_k$ for $H$ so that
each $\p D^{\d}_k$ intersects the dividing set $\G_{\p H}$ of $\p
H$ at two points and $H \setminus \cup_k D^{\d}_k$ is a standard
contact $3$-ball, after rounding the edges.
\end{itemize}
\end{itemize}}

\vspace{0.1in}

{\Def \label{stball} A standard contact $3$-ball is a tight
contact $3$-ball with a convex boundary whose dividing set is
connected.}

\vspace{0.2in}

Based on Theorem~\ref{honkazmat}, Honda, Kazez and Mati\'{c}
describe a \pob on $(M,\G)$ in Section 2 of their article
\cite{hkm1}. In this paper, for the sake of simplicity and without
loss of generality, we will assume that $M$ is connected. As a
consequence $M\setminus N(K)$ in Theorem~\ref{honkazmat} is also
connected.

In  \cite{eo}, we showed that the Honda-Kazez-Mati\'{c}
description gives an abstract \pob $(S,P,h)$, the balanced sutured
manifold associated to $(S,P,h)$ is isotopic to $(M,\G)$, and
$\xi$ is compatible with $(S,P,h)$. Conversely, let $(S,P,h)$ be
an abstract partial open book decomposition, $(M,\G)$ be the
balanced sutured manifold associated to it, and $\xi$ be a
compatible contact structure. Then we showed \cite{eo} that
$(S,P,h)$ is given by the Honda-Kazez-Mati\'c description.

\section{The EH-contact class is combinatorial}

The main result of \cite{hkm1} is the following:

{\Thm [\cite{hkm1}, Theorem 0.1]  Let $(M,\G)$ be a balanced
sutured manifold and let $\xi$ be a contact structure on $M$ with
convex boundary whose dividing set  on $\p M$ is isotopic to $\G$.
Then there exists an invariant $EH(M,\G,\xi)$ of the contact
structure $\xi$ which lives in $\SF/\{\pm1\}$. }

{\Rem Given a balanced sutured manifold $(M,\G)$, there exists a
contact structure $\xi$ on $M$ which makes $\p M$ convex and
realizes $\G$ as its diving set on $\p M$. Conversely given a
contact $3$-manifold $(M, \xi)$ (with convex boundary) whose
diving set is denoted by $\G$ on $\p M$, then $(M, \G)$ is a
balanced sutured manifold.}

{\Rem The $\pm 1$ ambiguity in the definition of $EH(M,\G,\xi)$ is
usually suppressed. An alternative way is to work with $\bfz_2$
coefficients.}

\vspace{0.1in}

Given a \pob $(S,P,h)$ consider the associated balanced sutured
manifold $(M, \G)$ and the uniquely (up to isotopy) determined
compatible contact structure $\xi$  on $M$. In this section we will
provide an algorithm to calculate the sutured Floer homology
$SF(-M,-\G)$ and the contact class $EH(M, \G,
\xi)$ in $\SF$ starting from $(S,P,h)$.

We now review basic definitions and properties of Heegaard
diagrams of sutured manifolds (cf. \cite{juh}).  A sutured
Heegaard diagram is given by $(\S, \A, \B)$, where the Heegaard
surface $\S$ is a compact oriented surface with nonempty boundary
and
 $\A= \{\a_1, \a_2, \ldots, \a_m \} $ and $\B = \{\b_1, \b_2, \ldots, \b_n \}
$ are two sets of pairwise disjoint simple closed curves in $\S
\setminus \p \S$. Every sutured Heegaard diagram  $(\S, \A, \B)$,
uniquely defines a sutured manifold $(M, \G)$ as follows: Let $M$
be the $3$-manifold obtained from $\S \times [0,1]$ by attaching
$3$-dimensional $2$-handles along the curves $\a_i \times \{0\}$
and $\b_j \times \{1\}$ for $i=1, \ldots, m$ and $ j= 1, \ldots,
n$. The suture $\G$ on $\p M$ is defined to be the set of curves $\p
\S \times \{ 1/2 \}$.

In \cite{juh}, Juh\'asz proved that if  $(M, \G)$ is defined by
$(\S, \A, \B)$, then $(M, \G)$ is balanced if and only if
$|\A|=|\B|$, the surface $\S$ has no closed components and both
\Atek and \Btek consist of curves linearly independent in $H_1 (
\S , \mathbb{Q})$. Hence a sutured Heegaard diagram $(\S, \A, \B)$
is called balanced if it satisfies the conditions listed above. We
will abbreviate balanced sutured Heegaard diagram as balanced
diagram from now on.

A partial open book decomposition of $(M,\G)$ gives a sutured
Heegaard diagram $(\S,\A,\B)$ of $(M,-\G)$ as follows: Let
$$\Sigma=P \times \{ 0\}
 \cup - S \times \{-1\} / \sim \ \subset \p H$$
be the Heegaard surface. Observe that, modulo identifications,
$$\p \S =
(\overline{\p P \setminus \p S}) \times \{0\} \cup -(\overline{\p
S \setminus \p P}) \times \{-1\} \simeq -\G \ . $$ As in the proof
of Proposition~\ref{torisu}, let $a_1, a_2, \ldots , a_r$ be
properly embedded pairwise disjoint arcs in $P$ with endpoints on
$A$ such that $S \setminus \cup_j a_j$ deformation retracts onto
$\overline{S\setminus P}$. Then define two families $\A = \{
\alpha_1, \a_2, \ldots , \alpha_r \}$ and $\B = \{ \beta_1, \b_2,
\ldots , \beta_r\}$ of simple closed curves in the Heegaard
surface $\S$ by $\alpha_j = a_j \times \{ 0\} \cup a_j \times
\{-1\} / \sim $ and $\beta_j = b_j \times \{ 0\} \cup h(b_j)
\times \{ -1 \} / \sim $, where $b_j$ is an arc isotopic to $a_j$
by a small isotopy such that
\begin{itemize}
\item the endpoints of $a_j$ are isotoped along $\p S$, in the
direction given by the boundary orientation of $S$, \item $a_j$
and $b_j$ intersect transversely in one point $x_j$ in the
interior of $S$, \item if we orient $a_j$, and $b_j$ is given the
induced orientation from the isotopy, then the sign of the
intersection of $a_j$ and $b_j$ at $x_j$ is $+1$.
\end{itemize}
$(\Sigma, \A, \B)$ is a sutured Heegaard diagram of $(M,-\G)$.
Here the suture is $-\G$ since $\p \S$ is isotopic to $-\G$.

Next we would like to review the definition of the sutured Floer
homology $SFH(M, \G)$ given by Juh\'asz (for more details see
\cite{juh}). Let $(M, \G)$ be a balanced sutured manifold and
$(\S, \A, \B)$ be an admissible balanced diagram defining it. Then
$SFH(M, \G)$ is defined to be the homology of the chain complex
$(CF(\S , \A,\B),\p)$, where $CF(\S , \A,\B)$ is the free abelian
group generated by the points in
$$\mathbb{T}_{\a} \cap \mathbb{T}_{\b} = (\a_1 \times \a_2 \times \cdots \times
\a_r )\cap (\b_1 \times \b_2 \times \cdots \times \b_r ) \subset
Sym^r(\S).
$$
For ${\bf x}, {\bf y} \in \mathbb{T}_{\a} \cap \mathbb{T}_{\b}$,
let $\mathcal{M}_{{\bf x}, {\bf y}}$ denote
the moduli space of pseudo-holomorphic maps $$ u : \mathbb{D} = \{z
\in \bfc : |z| \leq 1\} \to Sym^r(\S)$$ satisfying

\begin{enumerate}
 \item $u(1) = \bf{x}$ and $u(-1) = \bf{y}$,
 \item $u(\p \mathbb{D} \cap \{ z \in \bfc : \mbox{Im} z \geq 0\} )\subset \mathbb{T}_{\a} $
 and $u(\p \mathbb{D} \cap \{ z \in \bfc : \mbox{Im} z \leq 0\}) \subset \mathbb{T}_{\b} $,
 \item $u(\mathbb{D}) \cap (\p \S \times Sym^{r-1} (\S)) =
 \emptyset$.
\end{enumerate}

\noindent Then the boundary map $\p$ is defined by
$$\p {\bf x} = \sum_{ \mu( {\bf x}, {\bf y} )=1} \# (\mathcal{M}_{{\bf x}, {\bf y}}) {\bf y}
$$ where $\mu( {\bf x}, {\bf y} )$ is the relative Maslov index of the
pair and $\# (\mathcal{M}_{{\bf x}, {\bf y}})$ is a signed count
of points in the $0$-dimensional
quotient (by the natural $\mathbb{R}$-action) of
$\mathcal{M}_{{\bf x}, {\bf y}}$.

Let $(S,P,h)$ be a \pob and let  $(M, \G)$ be the associated
balanced sutured manifold. In Section~\ref{def} we described a balanced diagram $(\S, \A, \B)$
defining $(M, -\G)$. By changing the order of \Atek and \Btek we obtain a balanced diagram of $(-M,-\G)$.

The balanced diagram $(\Sigma, \B, \A)$ is shown to be admissible
in \cite{hkm1}. Hence the sutured Floer homology group $SFH(-M,
-\G)$ can be defined using this diagram. The contact class
$EH(M,\G,\xi)$ is defined \cite{hkm1} to be the homology class in
$\SF$ which descends from the cycle ${\bf x}$ in the complex
$CF(\S , \B,\A)$, where ${\bf x}=(x_1 , x_2, \ldots , x_r ) \in
Sym^r (\S)$.

{\Thm\label{comb} Let $(S,P,h)$ be a partial open book decomposition. Then the
sutured Floer homology group $\SF$ and the contact class
$EH (M, \G, \xi)$ in $\SF$ can be calculated
combinatorially, where $(M,\G)$ is the balanced sutured manifold
associated to  $(S,P,h)$ and $\xi$ is a contact structure on $M$
compatible with $(S,P,h)$. }
\begin{proof}
A balanced diagram $(\Sigma, \A, \B)$ is called \emph{simple} if
every component of $\S \setminus (\A \cup \B)$ whose closure is
disjoint from $\p \S$ is a bigon or a square. In \cite{juh1},
Juh\'azs proves, by modifying the procedure of Sarkar and Wang
\cite{sw}, that any balanced diagram can be turned into a simple
one using some isotopies and handle slides of the \Atek and \Btek
curves on $\S$. This provides the first step of an algorithm to
calculate the sutured Floer homology combinatorially since the
boundary homomorphism in the chain complex defining the homology
induced by a
simple balanced diagram can be calculated combinatorially. In
fact, exactly along the same lines as in the proof of Theorem 2.1
in \cite{p} one can see that, in our situation, i.e. when the
balanced diagram is obtained from a partial open book decomposition
as above, no
handle slide is necessary and the diagram can be modified into a
simple one by a sequence of isotopies on $P \times \{ 0\} \subset
\Sigma$ away from $A$. Denote the composition of these isotopies
by $\phi$ and observe that $\phi$ is a diffeomorphism fixing $A$
and isotopic to identity. The resulting simple diagram corresponds
to the partial open book decomposition $(S,P,h')$, where $h'=\phi \circ h$.
Hence the cycle ${\bf x} \in Sym^r(S\times \{ -1\}) \subset Sym^r
(\Sigma)$ considered in this simple balanced diagram still
descends to $EH(M,\G,\xi) $.

Once we have a simple diagram, by \cite{juh1}, it is combinatorial
to calculate the
boundary map of the sutured Floer chain complex. We just make a list
of all the
generators and count all the empty embedded bigons and squares on the
Heegaard surface connecting these generators by examining the diagram.
Finally by using simple linear algebra we can compute $\SF$ and identify
$[{\bf x}] = EH(M,\G,\xi) \in \SF$.
\end{proof}

In the rest of this section we present a few examples to demonstrate the
procedure explained in the proof of Theorem~\ref{comb} above.

\vspace{0.1in}

\begin{figure}[ht]
  \relabelbox \small {
  \centerline{\epsfbox{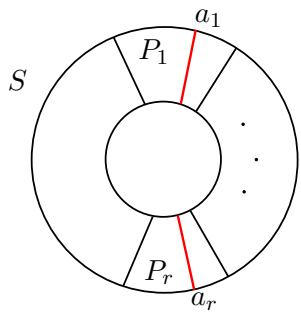}}}
  \relabel{1}{{$S$}}
  \relabel{2}{{$P_1$}}
 \relabel{3}{{$P_r$}}
 \relabel{4}{{$a_1$}}
  \relabel{5}{{$a_r$}}

  \endrelabelbox
        \caption{The annulus $S$, $r$ components $P_1, \dots , P_r$
        of $P$, and a basis $\{ a_1, \dots , a_r \}$ in Example 1.}
        \label{ex1-1}
\end{figure}

{\bf{Example 1.}} Let $S$ be an annulus, $P$ be a regular
neighborhood of $r$ disjoint and homotopically trivial arcs
connecting the two distinct boundary components of $S$, and the
monodromy $h$ be the inclusion of $P$ into $S$ (see
Figure~\ref{ex1-1}).

\begin{figure}[ht]
  \relabelbox \small {
  \centerline{\epsfbox{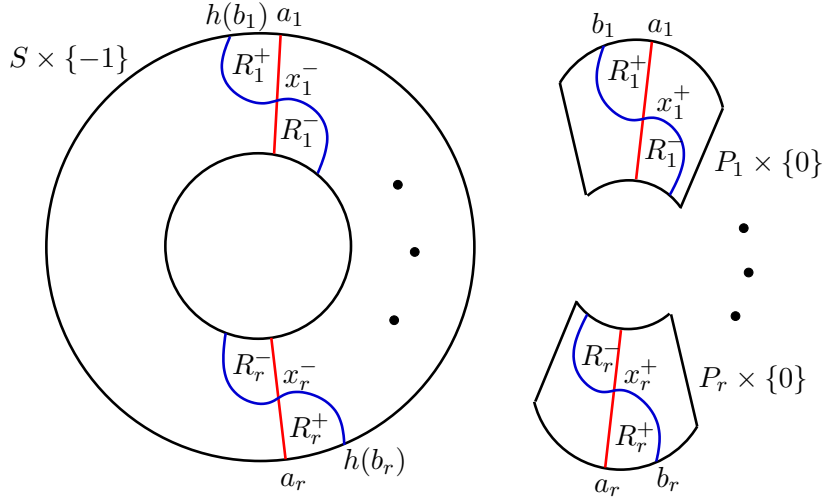}}}
  \relabel{1}{{$S \times \{-1\}$}}
  \relabel{2}{{$h(b_1)$}}
   \relabel{3}{{$a_1$}}
  \relabel{4}{{$h(b_r)$}}
   \relabel{5}{{$a_r$}}
   \relabel{6}{{$x_r^-$}}
  \relabel{7}{{$x_1^-$}}
   \relabel{8}{{$R_1^+$}}
  \relabel{9}{{$R_1^-$}}
   \relabel{10}{{$R_r^-$}}
   \relabel{11}{{$R_r^+$}}
  \relabel{12}{{$P_1 \times \{ 0\}$}}
   \relabel{13}{{$P_r \times \{ 0\}$}}
  \relabel{14}{{$a_1$}}
   \relabel{15}{{$b_1$}}
   \relabel{16}{{$a_r$}}
  \relabel{17}{{$b_r$}}
   \relabel{18}{{$x_1^+$}}
  \relabel{19}{{$x_r^+$}}
   \relabel{20}{{$R_1^+$}}
  \relabel{21}{{$R_1^-$}}
   \relabel{22}{{$R_r^-$}}
   \relabel{23}{{$R_r^+$}}
  \endrelabelbox
        \caption{ The surfaces $S \times \{-1\}$ and  $P \times \{0\}$,
        \Atek and \Btek curves, their intersections $x_j^{\pm}$, and the regions $R_j^{\pm}$ in Example 1.}
        \label{ex1-2}
\end{figure}

According to the notation in Figure~\ref{ex1-2}, the chain complex
$CF(\S , \B,\A)$ is generated by the $2^r$ generators $\{
(x_1^{\eps_1} , x_2^{\eps_2}, \dots, x_r^{\eps_n}) \}$, where
$\eps_j =\pm $ for $j=1,2, \ldots, r$.  On the other hand, the
regions that contribute to the boundary homomorphism $\p$ are
$R_1^{\pm},R_2^{\pm}, \ldots , R_r^{\pm}$. Each bigon $R^{\pm}_j$
effects only the generators of the form $(x_1^{\eps_1} ,
x_2^{\eps_2}, \dots, x_j^-, \dots , x_r^{\eps_n})$ and the
contribution is $\pm 1$ times the generator which differs only in
the $j^{th}$ component. The fact that the contribution has
absolute value $1$ follows from Theorem 7.4 in \cite{juh1} and for
each $j$ the signs of the contributions of $R^{\pm}_j$ are
opposite of each other by Lemma 9.1 (and especially the part of
its proof regarding the choice of a coherent system of
orientations) in \cite{os1}. For example, $\p (x_1^-, x_2^+,
x_3^+, \dots , x_r^+) = (x_1^+, x_2^+, x_3^+, \dots , x_r^+) -
(x_1^+, x_2^+, x_3^+, \dots , x_r^+)$, where the first term is
induced by $R_1^+$ and the second term is induced by $R_1^-$.
Consequently, the boundary map is trivial, hence $\SF \cong CF(\S
, \B,\A) \cong \bfz^{2^r}$, and $$ EH(M,\G,\xi) = [{\bf x}] =
[(x_1^+, x_2^+ , \dots ,x_r^+)]$$ is a generator of one of the
$\bfz$ summands.

Next we would like to describe the balanced sutured manifold
$(M,\G)$ associated to this partial open book decomposition
$(S,P,h)$, where we fix a positive integer $r$ for the rest of the
discussion. For each positive integer $m$, let $Y(m)$ denote the
balanced sutured manifold obtained by taking out $m$ disjoint open
$3$-balls from a closed $3$-manifold $Y$ and declaring the suture
to have exactly one connected component on each component of $\p
Y(m)$, as in \cite{juh}. Then we claim that $(M,\G) \cong Y(r)$
for $Y=S^1 \times S^2$. To prove our claim we observe that the
closed $3$-manifold which corresponds to the open book
decomposition with an annulus page and identity monodromy is $S^1
\times S^2$. Thus $M$ can be obtained from $Y$ by taking out $r$
disjoint open $3$-balls corresponding to $r$ connected components
of $S \setminus P$. Moreover by our construction in
Section~\ref{def} the suture $\G$ has $r$ connected components
each of which belongs to a different component of $\p M$. In the
light of this observation, the sutured Floer homology can be
calculated alternatively by using Proposition 9.14 in \cite{juh}
which states that $SFH(Y(m)) \cong \bigoplus_{2^{m-1}}
\widehat{HF} (Y)$ and the fact that $\widehat{HF}(S^1 \times
S^2)\cong \bfz \oplus \bfz$.

Furthermore we can identify the contact structure $\xi$ on $M$
which is compatible with $(S,P,h)$ as the contact structure
obtained by removing $r$ disjoint standard contact open $3$-balls
from the unique (up to isotopy) tight contact structure
$\xi_{std}$ on $S^1 \times S^2$. Hence the nontriviality of
$EH(M,\G,\xi)$ also follows from Theorem 4.5 in \cite{hkm1}
combined with the facts that $$EH(S^1 \times S^2 , \xi_{std}) =
c(S^1 \times S^2 , \xi_{std}) \;\; \mbox{and} \;\; c(S^1 \times
S^2 , \xi_{std}) \neq 0 , $$ where $c(S^1 \times S^2 , \xi_{std})$
denotes the contact Ozsv\'{a}th-Szab\'{o} invariant.

\vspace{0.1in}

\begin{figure}[ht]
  \relabelbox \small {
  \centerline{\epsfbox{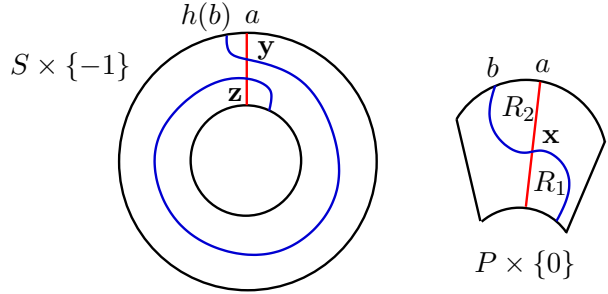}}}
  \relabel{1}{{$S \times \{-1\}$}}
  \relabel{2}{{$h(b)$}}
   \relabel{3}{{$a$}}
     \relabel{4}{{${\bf z}$}}
  \relabel{5}{{${\bf y}$}}
   \relabel{6}{{$R_2$}}
  \relabel{7}{{$R_1$}}
  \relabel{12}{{$P \times \{ 0\}$}}
     \relabel{11}{{$a$}}
   \relabel{10}{{$b$}}
  \relabel{9}{{${\bf x}$}}

  \endrelabelbox
        \caption{ The surfaces $S \times \{-1\}$ and $ P \times \{0\}$,
        \Atek and \Btek curves, their intersections ${\bf x},{\bf y},{\bf z}$, and the
        regions $R_1$ and $R_2$ in Example 2.}
        \label{ex2}
\end{figure}

{\bf{Example 2.}} Let $S$ and $P$ be as in the previous example
for $r=1$ and the monodromy $h$ be the restriction (to $P$) of a
{\emph{left}}-handed Dehn twist along the core of $S$. Using the
notation in Figure~\ref{ex2}, the generators of the chain complex
are ${\bf x},{\bf y}$ and ${\bf z}$. Moreover $\p {\bf x} =0$, $\p
{\bf y} = {\bf x}$ (by $R_1$) and $\p {\bf z}= {\bf x}$ (by
$R_2$). Hence $\SF = \bfz$ and $EH(M,\G,\xi)=0$. This is
consistent with the fact that the open book decomposition with
annulus page and left-handed Dehn twist monodromy is compatible
with an overtwisted contact $S^3$. In addition, by Proposition 4.2
in \cite{hkm1}, if the monodromy of a partial open book
decomposition is not {\emph{right-veering}}, i.e., if there is a
properly embedded arc $l \subset P$ with endpoints on $A$ such
that $h(l)$ is not to the right of $l$, then the contact invariant
of the compatible contact structure is zero.

\vspace{0.1in}

\begin{figure}[ht]
  \relabelbox \small {
  \centerline{\epsfbox{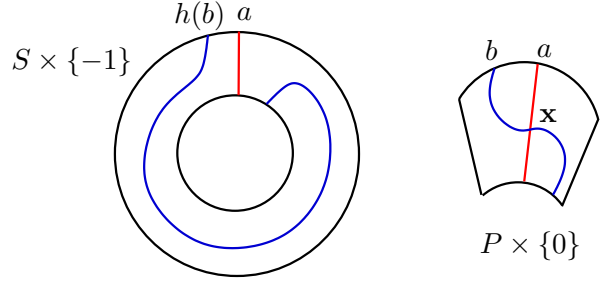}}}
  \relabel{1}{{$S \times \{-1\}$}}
  \relabel{2}{{$h(b)$}}
   \relabel{3}{{$a$}}
    \relabel{12}{{$P \times \{ 0\}$}}
     \relabel{11}{{$a$}}
   \relabel{10}{{$b$}}
  \relabel{9}{{${\bf x}$}}

  \endrelabelbox
        \caption{ The surfaces $S \times \{-1\}$ and $ P \times \{0\}$,
        \Atek and \Btek curves, and the intersection ${\bf x}$ for $r=1$ in Example 3.}
        \label{ex3}
\end{figure}

{\bf{Example 3.}} (\textbf{Standard contact $3$-ball}) Let $S$ and
$P$ be as in the first example for $r=1$, and the monodromy $h$ be
the restriction (to $P$) of a {\emph{right}}-handed Dehn twist
along the core of $S$. Then using the notation in
Figure~\ref{ex3}, there is a single generator ${\bf x}$ in the
chain complex $CF(\S , \B,\A)$ and hence the boundary homomorphism
is trivial. It follows that  $\SF \cong \bfz$ and $EH(M,\G,\xi)$
is a generator.

\begin{figure}[ht]
  \relabelbox \small {
  \centerline{\epsfbox{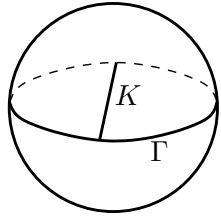}}}

  \relabel{2}{{$K$}}
 \relabel{3}{{$\G$}}

  \endrelabelbox
        \caption{The Legendrian arc $K$ in the standard contact $3$-ball.}
        \label{poball}
\end{figure}

\begin{figure}[ht]
  \relabelbox \small {
  \centerline{\epsfbox{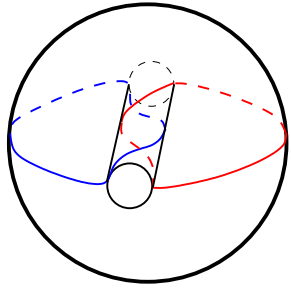}}}

  \endrelabelbox
        \caption{The dividing curves  on $\p H$.}
        \label{divcurve}
\end{figure}

In fact we can identify the contact $3$-manifold $(M, \G, \xi)$ as
the standard contact $3$-ball (cf. Definition~\ref{stball}). Here
the Legendrian graph $K$ which satisfies the conditions in
Theorem~\ref{honkazmat} is a single arc in $B^3$ connecting two
distinct points on $\G$ as depicted in Figure~\ref{poball}. The
complement $H$ of a regular neighborhood $N = N(K)$ in the
standard contact $3$-ball $B^3$ is a solid torus with two parallel
 dividing curves (see Figure~\ref{divcurve}) on $\p H$ which are homotopically
 nontrivial inside $H$. Here a meridional disk in $H$ will serve as the
 required compressing disk  $D_1^\d$ for $H$ in Theorem~\ref{honkazmat}
 $(B)$. On the other hand, $N$ is already a standard contact $3$-ball.

This shows in particular that the standard contact $3$-ball can be
obtained from a tight solid torus $H$ by attaching a tight
$2$-handle $N$. Note that, in Example 2, we attached a tight
$2$-handle to a tight solid torus to obtain an overtwisted contact
structure.

With a different point of view, one can see that $(M,\G) \cong
S^3(1)$, which is obtained by removing an open $3$-ball from the
standard tight contact $S^3$. Observe that the open book
decomposition with annulus page and right-handed Dehn twist
monodromy is compatible with the standard tight contact $S^3$ and
hence has nonzero contact class. Therefore Theorem 4.5 in
\cite{hkm1} already implies that the contact invariant in this
example is not zero. Moreover, the sutured Floer homology
calculations are consistent with the aforementioned result of
Juh\'asz.

\begin{figure}[ht]
  \relabelbox \small {
  \centerline{\epsfbox{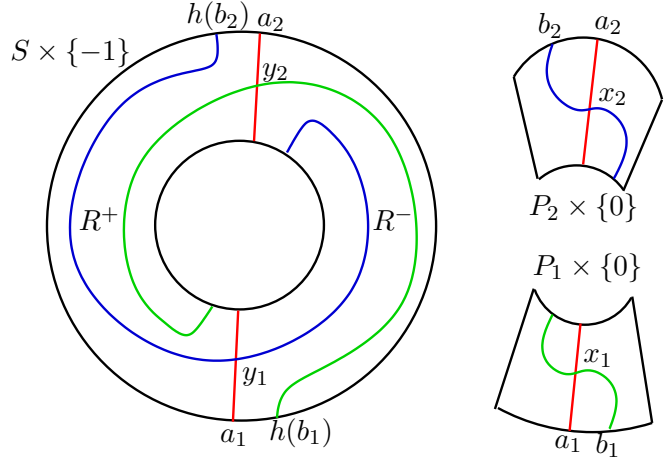}}}
  \relabel{1}{{$S \times \{-1\}$}}
  \relabel{2}{{$h(b_2)$}}
   \relabel{3}{{$a_2$}}
   \relabel{5}{{$h(b_1)$}}
   \relabel{4}{{$a_1$}}
    \relabel{12}{{$P_2 \times \{ 0\}$}}
     \relabel{6}{{$y_2$}}
   \relabel{10}{{$b_2$}}
  \relabel{7}{{$y_1$}}
  \relabel{17}{{$R^+$}}
   \relabel{18}{{$R^-$}}
\relabel{11}{{$a_2$}}
   \relabel{14}{{$a_1$}}
   \relabel{15}{{$b_1$}}
   \relabel{9}{{$x_2$}}
   \relabel{16}{{$x_1$}}
   \relabel{13}{{$P_1 \times \{0\}$}}

  \endrelabelbox

        \caption{ The surfaces $S \times \{-1\}$ and $ P \times \{0\}$,
        \Atek and \Btek curves, the intersections {\bf x} and {\bf y}, and the regions $R^{\pm}$
         for $r=2$ in Example 4.}
        \label{ex3-2}
\end{figure}

{\bf{Example 4.}} Let $S$ and $P$ be as in the first example for
$r \in \{2,3 \}$, and the monodromy $h$ be the restriction (to
$P$) of a {\emph{right}}-handed Dehn twist along the core of $S$.
First consider the case $r=2$. Using the notation in
Figure~\ref{ex3-2}, the generators of the chain complex are ${\bf
x}=(x_1,x_2)$ and ${\bf y}=(y_1,y_2)$ with $\p {\bf x} =0$ and $\p
{\bf y} = {\bf x} - {\bf x}=0 $ (by $R^{\pm}$), where the opposite
signs for the contributions of $R^{\pm}$ follow from Lemma 9.1 in
\cite{os1} as in Example~1. Hence $\SF \cong \bfz \oplus \bfz$ and
$EH(M,\G,\xi)$ is a generator of one of the $\bfz$ summands. Note
that $(M,\G) \cong S^3(2)$.

\begin{figure}[ht]
  \relabelbox \small {
  \centerline{\epsfbox{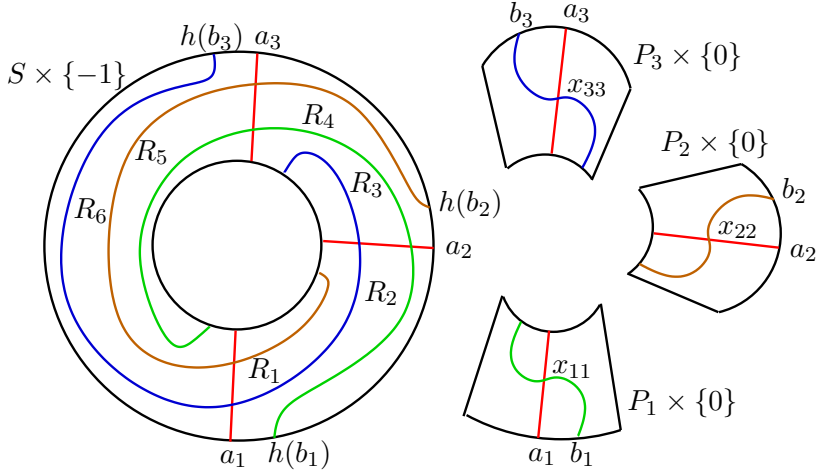}}}
  \relabel{1}{{$S \times \{-1\}$}}
  \relabel{2}{{$h(b_3)$}}
   \relabel{3}{{$a_3$}}
     \relabel{18}{{$h(b_2)$}}
   \relabel{17}{{$a_2$}}
   \relabel{5}{{$h(b_1)$}}
   \relabel{4}{{$a_1$}}
    \relabel{12}{{$P_2 \times \{ 0\}$}}
     \relabel{22}{{$P_3 \times \{ 0\}$}}
   \relabel{10}{{$b_2$}}
    \relabel{20}{{$b_3$}}
   \relabel{11}{{$a_2$}}
   \relabel{21}{{$a_3$}}
   \relabel{14}{{$a_1$}}
   \relabel{15}{{$b_1$}}
   \relabel{9}{{$x_{22}$}}
    \relabel{19}{{$x_{33}$}}
   \relabel{16}{{$x_{11}$}}
   \relabel{13}{{$P_1 \times \{0\}$}}
    \relabel{23}{{$R_1$}}
   \relabel{24}{{$R_2$}}
   \relabel{25}{{$R_3$}}
   \relabel{26}{{$R_4$}}
   \relabel{27}{{$R_5$}}
   \relabel{28}{{$R_6$}}

  \endrelabelbox

        \caption{ The surfaces $S \times \{-1\}$ and $ P \times \{0\}$,
        \Atek and \Btek curves, the intersections ${\bf x}_{ijk}$, and the regions $R_i$
         for $r=3$ in Example 4.}
        \label{ex3-3}
\end{figure}

Next consider the case $r=3$. Using the notation in
Figure~\ref{ex3-3}, the six generators of the chain complex $CF(\S
, \B,\A)$ are $\{ {\bf x}_{ijk} = (x_{1i} , x_{2j} , x_{3k}) : \{
i,j,k \} = \{ 1,2,3 \} \}$, where $x_{ij}$ is the single
intersection point in $\a_i \cap \b_j$, and the contact class
${\bf x}$ is ${\bf x}_{123}=(x_{11},x_{22},x_{33}) \in Sym^3(\S)$.
The boundary homomorphism is given by $\p {\bf x}_{123} = 0$, $\p
{\bf x}_{213} = {\bf x}_{123}-{\bf x}_{123}=0$ (by $R_1 \cup R_2$
and $R_4 \cup R_5$), $\p {\bf x}_{321} = {\bf x}_{123}-{\bf
x}_{123}=0$ (by $R_5 \cup R_6$ and $R_2 \cup R_3$ ), $\p {\bf
x}_{132} = {\bf x}_{123}-{\bf x}_{123}=0$ (by $R_3 \cup R_4$ and
$R_1 \cup R_6$), $\p {\bf x}_{231} = {\bf x}_{321} + {\bf x}_{213}
+ {\bf x}_{132}$ (by $R_1$, $R_3$ and $R_5$), and $\p {\bf
x}_{312} = {\bf x}_{321} + {\bf x}_{213} + {\bf x}_{132}$ (by
$R_2$, $R_4$ and $R_6$). As a result $\SF \cong \bfz \oplus \bfz
\oplus \bfz \oplus \bfz$ and $EH(M,\G,\xi)$ is a generator of one
of the $\bfz$ summands. Note that $(M,\G) \cong S^3(3)$.

\vspace{0.1in}

\begin{figure}[ht]
  \relabelbox \small {
  \centerline{\epsfbox{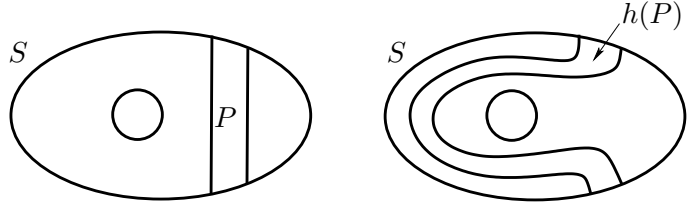}}}

  \relabel{2}{$h(P)$}

   \relabel{4}{$S$}

   \relabel{1}{$P$}
   \relabel{3}{$S$}
\endrelabelbox

        \caption{The partial open book decomposition $(S,P,h)$ in Example 5.}
        \label{ex4-1}
\end{figure}

{\bf{Example 5.}} (\textbf{Standard neighborhood of an overtwisted
disk}) Let $(S,P,h)$ be the partial open book decomposition shown
in Figure~\ref{ex4-1}.

\begin{figure}[ht]
  \relabelbox \small {
  \centerline{\epsfbox{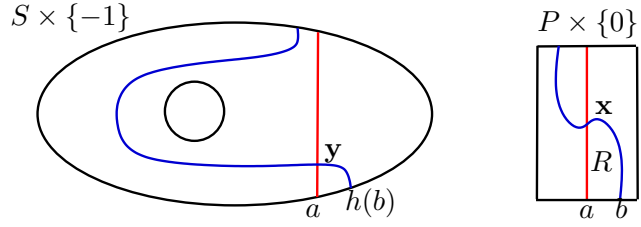}}}
  \relabel{7}{{$S \times \{-1\}$}}
  \relabel{2}{{${\bf y}$}}
   \relabel{8}{{$a$}}
   \relabel{9}{{$h(b)$}}
   \relabel{4}{{$a$}}
    \relabel{6}{{$P \times \{ 0\}$}}
     \relabel{5}{{$b$}}
   \relabel{1}{{${\bf x}$}}
   \relabel{3}{{$R$}}
\endrelabelbox

        \caption{ The surfaces $S \times \{-1\}$ and $ P \times \{0\}$,
        \Atek and \Btek curves, the intersections ${\bf x}$ and ${\bf y}$, and the region $R$
         in Example 5.}
        \label{ex4-2}
\end{figure}

Then using the notation in Figure~\ref{ex4-2}, the chain complex
$CF(\S , \B,\A)$ has two generators ${\bf x}$ and ${\bf y}$, and
the boundary homomorphism is given by $\p {\bf x} =0$, and $\p
{\bf y}={\bf x}$ by the bigon $R$. Hence $\SF = 0$ and obviously
$EH(M,\G,\xi) =0$. In fact, this is the partial open book
considered in Example 1 of \cite{hkm1} which is compatible with
the standard neighborhood of an overtwisted disk.

Here we observe that by Proposition~\ref{torisu}, $(M, \G, \xi)$
is obtained by gluing two compact connected contact $3$-manifolds
with convex boundaries, namely $(H, \G_{\p H} , \xi\vert_H)$ and
$(N, \G_{\p N} , \xi\vert_N)$ along parts of their boundaries. In
the following we will compute the sutured Floer homology groups
and the $EH$-classes of these contact submanifolds.

\begin{figure}[ht]
  \relabelbox \small {
  \centerline{\epsfbox{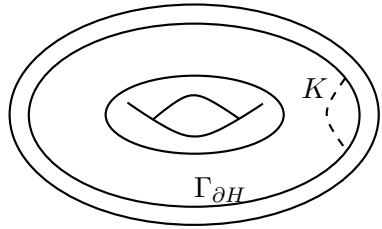}}}

   \relabel{1}{$K$}
 \relabel{3}{$\G_{\p H}$}
\endrelabelbox

        \caption{A Legendrian arc $K$ in $H$ connecting two points on $ \G_{\p H}$. }
        \label{solid}
\end{figure}

We know that
$$H=(S \times [-1,0])/\sim $$ where $S$ is an annulus and $(x,t)
\sim (x,t')$ for $x \in \p S$ and $t, t' \in [-1,0]$. There is a
unique (up to isotopy) compatible tight contact structure on $H$
whose dividing set $\G_{\p H}$ on $\p H$  is $\p S \times \{0\}$
(cf. Proposition~\ref{torisu}). Hence $(H, \G_{\p H}, \xi\vert_H)$
is a solid torus carrying a tight contact structure where $\G_{\p
H}$ consists of two parallel curves on $\p H$ which are
homotopically nontrivial in $H$. We observe that when we cut $H$
along a compressing disk we get a standard contact $3$-ball $B^3$
with its connected dividing set $\G_{\p B^3}$ on its convex
boundary. Note that $\G_{\p B^3}$ is obtained by ``gluing" $\G_{\p
H}$ and the dividing set on the compressing disk. Let $K$ be a
properly embedded Legendrian arc (as depicted in
Figure~\ref{solid}) in $B^3$ connecting two points on $ \G_{\p
H}$.

\begin{figure}[ht]
  \relabelbox \small {
  \centerline{\epsfbox{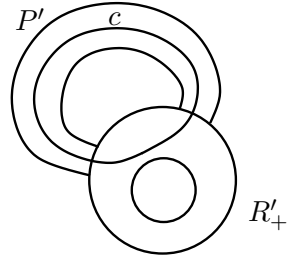}}}
 \relabel{2}{$R'_+$}
   \relabel{1}{$P'$}
   \relabel{3}{$c$}
\endrelabelbox

        \caption{The surface $S'= R'_+ \cup P'$, and the curve $c$ in $S'$.}
        \label{ex4-3}
\end{figure}

\begin{figure}[ht]
  \relabelbox \small {
  \centerline{\epsfbox{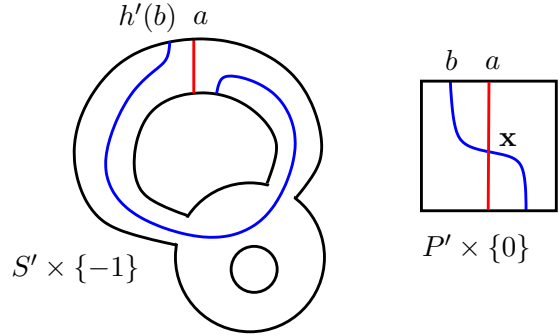}}}

  \relabel{2}{$P'\times \{0\}$}
   \relabel{1}{$S' \times \{-1\}$}
   \relabel{3}{$a$}
\relabel{4}{$h'(b)$}
 \relabel{5}{$a$}
\relabel{6}{$b$}
 \relabel{7}{{${\bf x}$}}
\endrelabelbox

        \caption{ The surfaces $S' \times \{-1\}$ and $ P' \times \{0\}$,
        \Atek and \Btek curves, the intersection ${\bf x}$.}
        \label{ex4-4}
\end{figure}

Then $K$ can be viewed as a Legendrian arc (disjoint from the
compressing disk) in $H$ connecting two distinct points of $\G_{\p
H}$ which satisfies the conditions in Theorem~\ref{honkazmat} just
as we discussed in Example 3. This gives a partial open book
decomposition $(S', P', h')$ compatible with $(H, \G_{\p H},
\xi\vert_H)$. The page $S'$ is a thrice punctured sphere that can
be obtained by attaching a $1$-handle $P'$ to the annulus $R'_+ =
S$. The monodromy map $h'$ is the restriction onto $P'$ of a
right-handed Dehn twist along the curve $c$ which is shown in
Figure~\ref{ex4-3}. There is a single generator ${\bf x}$ as shown
in Figure~\ref{ex4-4} in the chain complex $CF(\S , \B,\A)$ and
thus $EH(M,\G,\xi)$ is a generator of $SFH (-H, -\G_{\p H}) \cong
\mathbb{Z}$.

Similarly we know that $$N=(P \times [0,1])/\sim$$ where $(x,t)
\sim (x,t')$ for $x \in A$ and $t, t' \in [0,1]$. There is a
unique (up to isotopy) compatible tight contact structure on $N$
whose dividing set $\G_{\p N}$ on $\p N$  is $\p P \times \{1/2\}$
(cf. Proposition~\ref{torisu}). We observe that $(N, \G_{\p N},
\xi\vert_N)$ is the standard contact $3$-ball. It follows that
$EH(N, \G_{\p N}, \xi\vert_N)$ is a generator of the sutured Floer
homology group $SFH (-N, -\G_{\p N}) \cong \mathbb{Z}$ (see
Example $3$).

\vspace{0.1in}

\begin{figure}[ht]
  \relabelbox \small {
  \centerline{\epsfbox{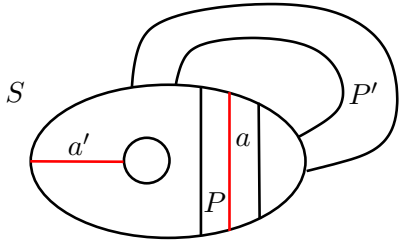}}}

\relabel{4}{$a$} \relabel{5}{$a'$}
 \relabel{1}{$P$}
   \relabel{2}{$P'$}
  \relabel{3}{$S$}
\endrelabelbox

        \caption{The page $T= S \cup P'$, and a basis $\{a, a'\}$ of $T$}
        \label{doubleo}
\end{figure}

{\bf{Example 6.}} Let us denote the neighborhood of an overtwisted
disk in Example 5 as  $(B^3 , \xi) $. We can glue two copies of
this overtwisted $3$-ball to get an overtwisted $(S^3, \xi_{ot})$.
In this example we will explicitly construct (cf. Theorem 4.5 in
\cite{hkm1}) the corresponding open book decomposition of $(S^3,
\xi_{ot})$ using the partial open book decomposition of $(B^3 ,
\xi) $. The page $T= S \cup P'$ is a twice punctured disk as
depicted in Figure~\ref{doubleo}. The monodromy is the composition
of a right-handed Dehn twist along one of the punctures and a
left-handed Dehn twist along the boundary of the disk. A simple
computation (see for example \cite{ozst}) shows that $d_3
(\xi_{ot}) = 1/2 $, which helps us identify the homotopy class of
$\xi_{ot}$ in $S^3$. A basis $ \{a, a'\}$ to compute $EH(S^3,
\xi_{ot})$ is also shown in Figure~\ref{doubleo}.

\begin{figure}[ht]
  \relabelbox \small {
  \centerline{\epsfbox{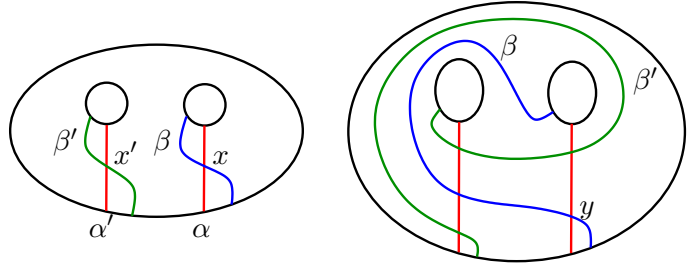}}}

\relabel{4}{$\a$}
 \relabel{1}{$x'$}
   \relabel{2}{$x$}
  \relabel{3}{$\a'$}
\relabel{5}{$\b'$} \relabel{6}{$\b$} 
\relabel{7}{$y$}  \relabel{c}{$\b'$}  \relabel{d}{$\b$}

\endrelabelbox

        \caption{A Heegaard diagram $(\S, \{\b,
\b'\}, \{\a, \a'\}) $ for the open book decomposition of $(S^3,
\xi_{ot})$. }
        \label{eh}
\end{figure}

In Figure~\ref{eh} we depicted a Heegaard diagram $(\S, \{\b,
\b'\}, \{\a, \a'\}) $ for the open book decomposition of $(S^3,
\xi_{ot})$. There are seven generators of $CF(\S , \B,\A)$ and the boundary relation $\p (y,x') = (x,x')$ implies that 
$EH(S^3, \xi_{ot})= [(x, x')]=0$, as expected.

Note that when we erase $\a', \b'$ and the handle $P'$ in
Figure~\ref{eh} we get a copy of Figure~\ref{ex4-2}.

\section{Gluing simple pieces}

In this section we generalize the discussion at the end of Example
5 to an arbitrary contact (connected) $3$-manifold $(M, \G, \xi)$
with the dividing set $\G$ on its convex boundary. It is clear
that $(M, \G, \xi)$ can be obtained by gluing two contact
handlebodies $(H, \G_{\p H}, \xi \vert_H)$  and $(N, \G_{\p N},
\xi\vert_N)$ along parts of their boundaries. We know that these
handlebodies are tight by Proposition~\ref{torisu} but we also
know that not every tight contact $3$-manifold has nontrivial
$EH$-class (cf. \cite{hkm1}). Nevertheless we have

{\Lem \label{nontrivial}  $EH (H, \G_{\p H}, \xi \vert_H)  \neq 0
$ and $EH (N, \G_{\p N}, \xi\vert_N) \neq 0 $}.

\begin{proof} We claim that $(H, \G_{\p H}, \xi\vert_H)$ can be embedded into
a Stein fillable \emph{closed} contact $3$-manifold $(Y, \xi')$.
We just embed $H$ into an open book decomposition (in the usual
sense) with page $S$ and trivial monodromy whose compatible
contact structure is Stein fillable by \cite{g} (and hence tight
by \cite{eg}). To be more precise, we embed $$H=(S \times
[-1,0])/\sim $$ where $(x,t) \sim (x,t')$ for $x \in \p S$ and $t,
t' \in [-1,0]$ into
$$Y = (S \times [-2,0])/\sim
$$ where $(x,0) \sim (x,-2)$ for $x \in S$ and $(x,t) \sim (x,t')$
for $x \in \p S$ and $t, t' \in [-2,0]$. Let $\xi'$ be the tight
structure on $Y$ which is compatible with the above open book
decomposition.  Then $\p H = S \times \{0 \} \cup - S \times
\{-1\}$ (which is obtained by gluing two pages along the binding)
can be made convex with respect to $\xi'$ so that the dividing set
on $\p H $ is exactly the binding (see \cite{e} for example).
Since the contact Ozsv\'{a}th-Szab\'{o} invariant $ c(Y, \xi')$ of
a Stein fillable contact $3$-manifold is nontrivial by \cite{os}
and $c(Y, \xi') = EH(Y, \xi')$ by \cite{hkm}, it follows by
Theorem 4.5 in \cite{hkm1} that $EH (H, \G_{\p H}, \xi\vert_H)
\neq 0 $. A similar argument shows that $EH (N, \G_{\p N},
\xi\vert_N) \neq 0 $. Indeed this is simply because $(N, \G_{\p
N}, \xi\vert_N)$ can be embedded in $(H, \G_{\p H}, \xi \vert_H)$
(see \cite{eo} for details).
\end{proof}

{\Cor \label{glu} Every contact $3$-manifold can be obtained by
gluing tight contact handlebodies (with convex boundaries) whose
$EH$-classes are nontrivial. }

\begin{proof} Suppose that $(Y, \xi)$ is a closed connected contact $3$-manifold.
We know that $(Y, \xi)$ has a compatible open book decomposition
with page $S$ and monodromy $h: S \to S$ (cf. \cite{g}). This
implies that $$ Y \cong (S \times [-2,0])/\sim $$ where $(h(x),0)
\sim (x,-2)$ for $x \in S$ and $(x,t) \sim (x,t')$ for $x \in \p
S$ and $t, t' \in [-2,0]$. Thus $(Y, \xi)$ is obtained by gluing
the contact handlebodies $(H_1, \G_{\p H_1}, \xi\vert_{H_1})$ and
$(H_2, \G_{\p H_2}, \xi\vert_{H_2})$ along their convex boundaries
using the monodromy map. Here $$H_1 \cong (S \times [-1,0])/\sim
$$ where $(x, 0) \sim (x,-1)$ for $x \in S$ and $(x,t) \sim
(x,t')$ for $x \in \p S$ and $t, t' \in [-1,0]$ and $$H_2 \cong (S
\times [-2,-1])/\sim $$ where $(x, -1) \sim (x,-2)$ for $x \in S$
and $(x,t) \sim (x,t')$ for $x \in \p S$ and $t, t' \in [-2,-1]$.
Then Lemma~\ref{nontrivial} implies that $EH (H_i, \G_{\p H_i},
\xi\vert_{H_i}) \neq 0 $ for $i=1,2$. Now suppose that $(Y, \xi)$
is connected contact $3$-manifold with nonempty convex boundary.
Then $Y$ admits a partial open book decomposition $Y\cong H \cup
N$ and the result is immediate by Lemma~\ref{nontrivial}.

If $Y$ has more than one connected components then we can apply
the above argument to each of its components to obtain the desired
result.
\end{proof}

{\Rem It is well-known that an overtwisted closed contact
$3$-manifold has trivial Ozsv\'{a}th-Szab\'{o} invariant
\cite{os}. Every such manifold is obtained by gluing two tight
contact handlebodies (with convex boundaries) whose $EH$-classes
are nontrivial by Corollary~\ref{glu}. }

\vspace{0.2in}

In Lemma~\ref{nontrivial} we proved that  $EH (H, \G_{\p H}, \xi
\vert_H)  \neq 0 $. In fact we have

{\Prop  \label{z} The contact class $EH (H, \G_{\p H},
\xi\vert_H)$ is a generator of the sutured Floer homology group
$SFH (-H, -\G_{\p H}) \cong \mathbb{Z}.$ }

\begin{proof}  Let $S$ be a connected genus $g$ surface with $n$
boundary components. Then there is a unique (up to isotopy)
compatible tight contact structure on $H=(S \times [-1,0])/\sim $,
where $(x,t) \sim (x,t')$ for $x \in \p S$ and $t, t' \in [-1,0]$,
whose dividing set $\G_{\p H}$ on $\p H$ is $\p S \times \{0\}$
(cf. Proposition~\ref{torisu}). In the following we will describe
a partial open book decomposition $(S', P', h')$ for $(H, \G_{\p
H}, \xi\vert_H)$. We observe that when we cut the solid handlebody
$H$ along compressing disks we get a standard contact $3$-ball
$B^3$. Note that $\G_{\p B^3}$ is obtained by gluing $\G_{\p H}$
and the dividing sets on the compressing disks. Let $K$ be a
properly embedded Legendrian arc in $B^3$ connecting two points on
$ \G_{\p H}$. Then $K$ can be viewed as a Legendrian arc in $H$
connecting two distinct points of $\G_{\p H}$ and it satisfies the
conditions in Theorem~\ref{honkazmat} by our construction (see
Example 3). Then $S'$ will be obtained from $R'_+ = S$ by
attaching a $1$-handle $P'$. Let $c$ be the indicated curve in
Figure~\ref{surfP}. Then $h'$ is the restriction to $P'$ of a
right-handed Dehn twists along $c$. It follows (just as in Example
5, Figure~\ref{ex4-4}) that there is a unique intersection point
of the \Atek and \Btek curves and consequently, $EH(H, \G_{\p H},
\xi\vert_H)$ is a generator of $SFH (-H, -\G_{\p H}) \cong
\mathbb{Z}$.

\begin{figure}[ht]
  \relabelbox \small {
  \centerline{\epsfbox{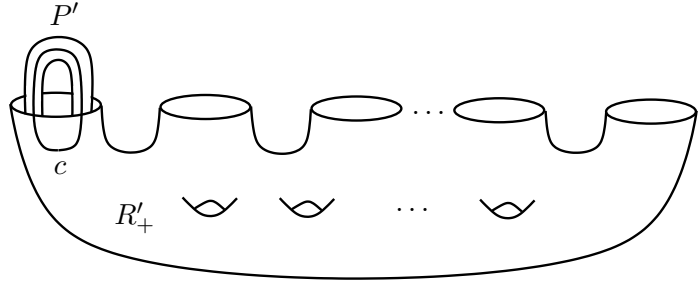}}}

\relabel{4}{$\dots$}
\relabel{5}{$\dots$}
 \relabel{1}{$P'$}
   \relabel{2}{$c$}
  \relabel{3}{$R'_+$}
\endrelabelbox

        \caption{The surface $S'= R'_+ \cup P'$, and the curve $c$ in $S'$.}
        \label{surfP}
\end{figure}
\end{proof}


The proof of the following result is very similar to the proof of
Proposition~\ref{z}.  Here we use the same notation as in
Lemma~\ref{nontrivial}.

{\Prop Suppose that $N_1, N_2, \dots , N_n$ are the connected
components of the handlebody $N$. Then for $1 \leq i \leq n$,
$$SFH (-N_i, -\G_{\p N_i}) \cong \mathbb{Z}.$$ Moreover $EH (N_i,
\G_{\p N_i}, \xi\vert_{N_i})$ is a generator of  $SFH (-N_i,
-\G_{\p N_i})$ and
$$EH (N, \G_{\p N}, \xi\vert_{N}) = \sum_{i=1}^{n} EH (N_i,
\G_{\p N_i}, \xi\vert_{N_i}) \in  \bigoplus_{i=1}^{n} SFH (-N_i,
-\G_{\p N_i}) \cong \mathbb{Z}^n.$$}

\vspace{0.2in}

\noindent{\bf {Acknowledgements.}} We would like to thank Andr\'as
Stipsicz and Sergey Finashin for valuable comments on a draft of
this paper. We would also like to thank John Etnyre for very
helpful email correspondence. TE was partially supported by a
GEBIP grant of the Turkish Academy of Sciences and a CAREER grant
of the Scientific and Technological Research Council of Turkey. BO
was partially supported by the research grant 107T053 of the
Scientific and Technological Research Council of Turkey.



\begin{thebibliography}{99999}
\bibitem{eg} Y. Eliashberg and M. Gromov,
{\em Convex symplectic manifolds,} Several complex variables and
complex geometry, Part 2 (Santa Cruz, CA, 1989), 135--162, Proc.
Sympos. Pure Math., 52, Part 2, Amer. Math. Soc., Providence, RI,
1991.


\bibitem{eo} T. Etg\"{u} and B. Ozbagci, {\em Relative Giroux
correspondence,} preprint, arXiv:math.GT/0802.0810

\bibitem{e}
J. B. Etnyre, {\em Lectures on open book decompositions and
contact structures,} Floer homology, gauge theory, and
low-dimensional topology, 103--141, Clay Math. Proc., 5, Amer.
Math. Soc., Providence, RI, 2006.



\bibitem{g} E. Giroux, \emph{G\'{e}ometrie de contact: de la dimension trois
vers les dimensions sup\'{e}rieures,} Proceedings of the
International Congress of Mathematicians (Beijing 2002), Vol. II,
405--414.

\bibitem{juh}
A. Juh\'{a}sz, {\em Holomorphic discs and sutured manifolds},
Algebr. Geom. Topol. 6 (2006), 1429--1457.

\bibitem{juh1} A. Juh\'{a}sz, {\em Floer homology and surface
decompositions}, preprint arXiv:math.GT/0609779



\bibitem{hkm}
K. Honda, W. Kazez and  G. Mati\'{c}, {\em On the contact class in
Heegaard Floer homology,} preprint, arXiv:math.GT/0609734

\bibitem{hkm1}
K. Honda, W. Kazez and  G. Mati\'{c}, {\em The contact invariant
in sutured Floer homology,} preprint, arXiv:math.GT/0705.2828v2

\bibitem{ozst}
B. Ozbagci and A. I. Stipsicz, {\em Surgery on contact
3--manifolds and Stein surfaces}, Bolyai Soc. Math. Stud., Vol.
{\bf 13}, Springer, 2004.


\bibitem{os1} P. Ozsv\'{a}th and Z. Szab\'{o}, {\em Holomorphic
disks and topological invariants for closed three-manifolds,} Ann.
of Math. (2) 159 (2004), no. 3, 1027--1158.

\bibitem{os} P. Ozsv\'{a}th and Z. Szab\'{o},
{\em Heegaard Floer homology and contact structures,} Duke Math.
J. 129 (2005), no. 1, 39--61.



\bibitem{p}
O. Plamenevskaya, {\em A combinatorial description of the Heegaard
Floer contact invariant}, Algebr. Geom. Topol. 7 (2007), 1201-1209.


\bibitem{sw}
S. Sarkar and J. Wang, {\em An algorithm for computing some
Heegaard Floer homologies}, preprint, arXiv:math.GT/0607777



\end{thebibliography}
\end{document}